\documentclass[12pt,reqno]{amsart}
\usepackage{amsmath, amssymb, amsthm}
\usepackage{url}
\usepackage{mathrsfs}
\usepackage[ansinew]{inputenc}
\usepackage[breaklinks]{hyperref}
\usepackage{fancyhdr}

\usepackage{multirow}

\allowdisplaybreaks

\renewcommand{\(}{\left\(}
\renewcommand{\)}{\right\)}
\renewcommand{\[}{\left\[}
\renewcommand{\]}{\right\]}
\newtheorem{remark}[]{Remark}
\numberwithin{equation}{section}
 \theoremstyle{plain}
\newtheorem{theorem}{Theorem}[section]

\newcommand{\sm}{\left(\begin{smallmatrix}}
\newcommand{\esm}{\end{smallmatrix}\right)}

   \makeatletter
\def\proof{\@ifnextchar[{\@oproof}{\@nproof}}
\def\@oproof[#1][#2]{\trivlist\item[\hskip\labelsep\textit{#2 \textbf{Proof of}\
#1.}~]\ignorespaces}
\def\@nproof{\trivlist\item[\hskip\labelsep\textit{Proof.}~]\ignorespaces}

\makeatother

\usepackage{color}
\usepackage{amsmath}

\definecolor{blue}{rgb}{0,0,1}
\definecolor{red}{rgb}{1,0,0}
\definecolor{green}{rgb}{0,.6,.2}
\definecolor{purple}{rgb}{1,0,1}

\long\def\red#1\endred{{\color{red}#1}}
\long\def\blue#1\endblue{{\color{blue}#1}}
\long\def\purple#1\endpurple{{\color{purple}#1}}
\long\def\green#1\endgreen{{\color{green}#1}}

\begin{document}
\title[On Euler's Theorem]{On Euler's Theorem}
\author{George E. Andrews}
\address{Department of Mathematics, The Pennsylvania State University, University Park, PA 16802, U.S.A.}
\email{gea1@psu.edu}

\author{Rahul Kumar}
\address{Department of Mathematics, Indian Institute of Technology, Roorkee-247667, Uttarakhand, India}
\email{rahul.kumar@ma.iitr.ac.in} 

\author{Ae Ja Yee}
\address{Department of Mathematics, The Pennsylvania State University, University Park, PA 16802, U.S.A.}
\email{yee@psu.edu}

\maketitle
\pagenumbering{arabic}
\pagestyle{headings}
\begin{abstract}
Euler's theorem asserts that $A(n)=B(n)$ where $A(n)$ is the number of partitions of $n$ into distinct parts and $B(n)$ is the number of partitions of $n$ into odd parts. In this paper, it is proved that for $n>0$,
\begin{align*}
A(n)=B(n)=C(n+1)=\frac{1}{2}D(n+1),
\end{align*}
where $C(n)$ is the number of partitions of $n$ with largest part even and parts not exceeding half of the largest part are distinct, and $D(n)$ is the number of partitions of $n$ into non-negative parts wherein the smallest part appear exactly twice and no other parts are repeated. 

\end{abstract}

\section{Introduction}

The grandfather of all partition identities is Euler's theorem \cite{andrews book}, namely,
\begin{align*}
A(n)=B(n),
\end{align*}
where $A(n)$ is the number of partitions of $n$ into distinct parts and $B(n)$ is the number of partitions of $n$ into odd parts.

This theorem contains all the elements that would suggest generalizations, and, over the centuries, generalizations have been found in profusion. The Rogers-Ramanujan identities \cite[p.~104]{andrews book} and Schur's 1926 theorem \cite[p.~116]{andrews book} kicked off the twentieth century's contributions, consult Henry Alder's survey \cite{alder} for an account of some of the results in the late 20th century. 

However, little has ever been added to Euler's theorem itself. In this paper, we shall add two further partition functions.

\begin{theorem}\label{main theorem}
For $n>0$, we have
\begin{align}\label{main theorem eqn}
A(n)=B(n)=C(n+1)=\frac{1}{2}D(n+1),
\end{align}
where $C(n)$ is the number of partitions of $n$ with largest part even and parts not exceeding half of the largest part are distinct, and $D(n)$ is the number of partitions of $n$ into non-negative parts wherein the smallest part appear exactly twice and no other parts are repeated.
\end{theorem}

For example, if $n=6$, the four sets of partitions are following:
\begin{center}
\begin{tabular}{ p{2cm} p{3.5cm} p{3cm} p{4.1cm} p{4.1cm}}

 \vspace{2mm}$A(6)$\newline  &  \vspace{2mm} $B(6)$ & \vspace{2mm} $C(7)$& \vspace{2mm} $D(7)$ \\

 $6$ \newline 5+1 \newline 4+2\newline 3+2+1   & 5+1 \newline 3+3 \newline 3+1+1+1 \newline 1+1+1+1+1+1  & 6+1\newline 4+3 \newline 4+2+1 \newline 2+2+2+1  & 0+0+7\newline 0+0+6+1\newline 0+0+5+2\newline 0+0+4+3\newline 0+0+4+2+1\newline 1+1+5\newline 1+1+2+3\newline 2+2+3
\end{tabular}
\end{center}

We conclude the introduction with the following remark.
\begin{remark}
In their recent paper \cite{ba}, M. El Bachraoui and the first author considered partitions with multiple appearances by the first part.  All parts were assumed to be positive. It would be a simple matter to extend the results of that paper to the case of non-negative parts in that this would add $(-q;q)_\infty$ to the generating functions in question.
\end{remark}

This paper is organised as follows. In section \ref{bc}, we provide the brief proof that $B(n)=C(n+1)$. In Section \ref{cd}, we prove that $C(n)=\frac{1}{2}D(n)$. We also provide bijective proofs our assertions in Theorem \ref{main theorem} in Section \ref{bijective proofs}.

\section{Proof that $B(n)=C(n+1)$}\label{bc}

We note that
\begin{align}\label{bn}
\sum_{n=1}^\infty B(n)q^{n+1}&=q\sum_{n=1}^\infty\frac{q^{2n-1}}{(q;q^2)_n} \notag \\
&=\sum_{n=1}^\infty\frac{q^{2n}}{(q;q^2)_n},
\end{align}
here, and throughout the sequel,
\begin{align*}
(A;q)_N=(1-A)(1-Aq)\cdots(1-Aq^{N-1}),
\end{align*}
and 
\begin{align*}
(A;q)_\infty=\lim_{N\to\infty}(A;q)_N.
\end{align*}

On the other hand, we have
\begin{align}
\sum_{n=1}^\infty C(n)q^{n}&=\sum_{n=1}^\infty\frac{(-q;q)_nq^{2n}}{(q^{n+1};q)_n}\nonumber\\
&=\sum_{n=1}^\infty\frac{(q^2;q^2)_nq^{2n}}{(q;q)_{2n}} \label{c1n}\\
&=\sum_{n=1}^\infty\frac{q^{2n}}{(q;q^2)_n}.\label{cn}
\end{align}

Comparing \eqref{bn} and \eqref{cn}, we see that $B(n)=C(n+1)$.

\section{Proof that $C(n)=\frac{1}{2}D(n)$}\label{cd}

Let us fix $C(0)=1$. From \eqref{c1n}, we have
\begin{align}\label{before p19}
\sum_{n=0}^\infty C(n)q^{n}&=\sum_{n=0}^\infty\frac{(q^2;q^2)_nq^{2n}}{(q;q)_{2n}}\nonumber\\
&=(q^2;q^2)_\infty\sum_{n=0}^\infty\frac{q^{2n}}{(q;q)_{2n}(q^{2n+2};q^2)_\infty}.
\end{align}
We now employ the following Euler's identity \cite[p.~19, (2.2.5)]{andrews book}
\begin{align}\label{euler identity}
\frac{1}{(t;q)_\infty}=\sum_{m=0}^\infty\frac{t^m}{(q,q)_m},
\end{align}
(with replacing $q$ by $q^2$ and then letting $t=q^{2n+2}$) in \eqref{before p19}  so as to obtain
\begin{align}\label{series}
\sum_{n=0}^\infty C(n)q^{n}&=(q^2;q^2)_\infty\sum_{n=0}^\infty\frac{q^{2n}}{(q;q)_{2n}}\sum_{m=0}^\infty\frac{q^{2nm+2m}}{(q^2;q^2)_m}\nonumber\\
&=(q^2;q^2)_\infty\sum_{m,n=0}^\infty\frac{q^{2n+2nm+2m}}{(q;q)_{2n}(q^2;q^2)_m}\nonumber\\
&=(q^2;q^2)_\infty\sum_{m,n=0}^\infty\frac{1}{2}\left(1+(-1)^n\right)\frac{q^{n+nm+2m}}{(q;q)_{n}(q^2;q^2)_m}\nonumber\\
&=\frac{1}{2}(q^2;q^2)_\infty\sum_{m=0}^\infty\frac{q^{2m}}{(q^2;q^2)_m}\left\{\sum_{n=0}^\infty\frac{q^{n(m+1)}}{(q;q)_n}+\sum_{n=0}^\infty\frac{(-1)^nq^{n(m+1)}}{(q;q)_n}\right\}.
\end{align}
Upon invoking \eqref{euler identity} twice, once with letting $t=q^{m+1}$ and once with letting $t=-q^{m+1}$, and then substituting both resulting expressions in \eqref{series}, we conclude that
\begin{align}\label{3.4}
\sum_{n=0}^\infty C(n)q^{n}&=\frac{1}{2}(q^2;q^2)_\infty\sum_{m=0}^\infty\frac{q^{2m}}{(q^2;q^2)_m}\left\{\frac{1}{(q^{m+1};q)_\infty}+\frac{1}{(-q^{m+1};q)_\infty}\right\}\nonumber\\
&=\frac{1}{2}(-q;q)_\infty\sum_{m=0}^\infty\frac{q^{2m}}{(-q;q)_m}+\frac{1}{2}(q;q)_\infty\sum_{m=0}^\infty\frac{q^{2m}}{(q;q)_m}\nonumber\\
&=\frac{1}{2}\sum_{m=0}^\infty q^{2m}(-q^{m+1};q)_\infty+\frac{1}{2}(1-q),
\end{align}
where the last step follows upon again using \eqref{euler identity} with letting $t=q^2$.

It is clear now that
\begin{align*}
\sum_{n=0}^\infty C(n)q^{n}&=\frac{1}{2}\sum_{n=0}^\infty D(n)q^{n}+\frac{1}{2}-\frac{q}{2}.
\end{align*}
Comparing the coefficients for $n>1$ in the above expression proves the last equality of \eqref{main theorem eqn}.

\bigskip

\section{A Bijective proof of Theorem \ref{main theorem}}\label{bijective proofs}

\subsection{A Bijective proof of $A(n-1)=\frac{1}{2}D(n)$}
We can first remove the non-negative part condition for $D(n)$ as follows:
\begin{quote}
For $n\ge 2$, $D(n)$ counts the number of partitions of $n$ where only the smallest part can repeat at most twice and all other parts are distinct. 
\end{quote}
In other words, a partition $\lambda=(\lambda_1,\ldots, \lambda_{\ell})$ counted by $D(n)$ satisfies
$$
\lambda_1>\lambda_2> \cdots > \lambda_{\ell-1}\ge \lambda_{\ell} \text{ if $\ell>1$}  
$$
or 
$\lambda=(\lambda_1)$ has only one part. 

We subtract $1$ from the smallest part $\lambda_{\ell}$ and call the resulting partition $\mu$. Then, it is clear that $\mu$ is a partition of $n-1$. Also, 
\begin{enumerate}
\item[] Case 1: if $\ell=1$, then $\mu$ has only one part; 

\item[] Case 2: if $\ell>1$ and $\lambda_{\ell-1}\ge \lambda_\ell>1$, then  $\mu$ has  $\ell$ distinct parts;

\item[] Case 3: if $\ell>1$ and $\lambda_{\ell-1}\ge \lambda_\ell=1$, then $\mu$ has $\ell-1$   distinct parts.
\end{enumerate}
The union of Cases 1 and 2 is the set of partitions of $n-1$ into distinct parts. Also, there are partitions of $n-1$ into distinct parts in Case 3. Therefore, we get $D(n) =2A(n-1)$ as desired.

\subsection{A Bijective proof of $B(n)=C(n+1)$}
This theorem follows from Glaisher's bijection $\phi$ for $A(n)=B(n)$ \cite{andrews, glaisher}. For completeness, we give a brief sketch of the bijection. For any positive integer $M$, we can write it as $M=2^k a$ for some $k\ge 0$ and some odd integer $a$. Then $\phi$ maps $2^k a$ to $2^k$ copies of $a$.  In other words, $\phi$ takes $N$ and splits it into halves repeatedly until there are no even parts left.

Let $\lambda$ be a partition counted by $C(n+1)$ with largest part equal to $2N$. We subtract $1$ from the largest part. So, the weight of $\lambda$ becomes $n$. 

Next, odd parts are allowed in parts counted by $B(n)$, so we need to take care of even parts of $\lambda$. Also, odd parts greater than $N$ can repeat while odd parts less than or equal to $N$ can appear at most once, so we will show that the application of Glaisher's bijection $\phi$ to even parts of $\lambda$ produces odd parts less than or equal to $N$ with multiplicity greater than $1$.  

Let $M \le 2N$ be an even number, which can be  written as
$$
M=2^k a  \text{ for some  odd integer $a$ and some  $ k\ge 1$}.
$$
We apply Glaisher's bijection $\phi$ to $M$ and obtain $2^k$ parts of size $a$. Note that since
$
2^k a \le 2N,
$
$$
a \le N.
$$
Also, if $M\le N$, then $M$ can appear only once, so we get all distinct positive powers $2^k$ such that
\begin{equation} 
2^k\le N/a. \label{power1}
\end{equation}
On the other hand, if $N< M\le 2N$, then
\begin{equation} 
2^k> N/a,  \label{power2}
\end{equation}
and $M$ can repeat. 
Suppose $M$ appears $f$ many times. Upon applying Glaisher's bijection $\phi$ to $f$ copies of $M$, we obtain
$$
 f 2^{k} \text{ copies of $a$}. 
$$
By writing $f$ as a binary expansion, 
$$
f 2^k =(f_0 \,2^{0}+f_1 \, 2^{1} +\cdots  ) 2^k, 
$$
where $f_j$ is either $0$ or $1$  for $j\ge 0$. By \eqref{power2}, we see that each summand in the above expression represents a distinct power of $2$ greater than $N/a$, i.e., 
\begin{equation}
f_j \, 2^{j+k} >N/a. \label{power3}
\end{equation}
It follows from \eqref{power1} and \eqref{power3} that $a$ can appear with any multiplicity greater than $1$.  This proves that the resulting partition is counted by $B(n)$. 



\medskip

\noindent {\bf{Acknowledgements:}}\, The first author is partially supported by the Simons Foundation Grant 633284, and the second author is partially supported by the Grant number ANRF/ECRG/2024/003222/PMS of Anusandhan National Research Foundation (ANRF), Govt. of India, and the FIG grants of IIT Roorkee. Both authors sincerely thank these institutions for their support.

\end{document}